\documentclass[11pt]{article}
\usepackage{amscd}
\usepackage{amsfonts}
\usepackage{amsmath}
\usepackage{amssymb}
\usepackage{amsthm}
\usepackage{bbm}
\usepackage{CJK}
\usepackage{fancyhdr}
\usepackage{graphicx}
\usepackage{hyperref}
\usepackage{indentfirst}
\usepackage{latexsym}
\usepackage{mathrsfs}
\usepackage{xypic}
\usepackage[compress]{cite}

\newtheorem{theorem}{Theorem}[section]
\newtheorem{lemma}[theorem]{Lemma}
\newtheorem{definition}[theorem]{Definition}
\newtheorem{proposition}[theorem]{Proposition}
\newtheorem{example}[theorem]{Example}

\newtheorem{remark}[theorem]{Remark}

\usepackage[top=1in,bottom=1in,left=1.25in,right=1.25in]{geometry}
\def\<{\langle}
\def\>{\rangle}

\def\o{\otimes}

\date{\today}
\begin{document}
\renewcommand{\baselinestretch}{1.2}
\renewcommand{\arraystretch}{1.0}
\title{\bf Deformations and  extensions of  modified $\lambda$-differential  Lie-Yamaguti algebras}
\date{\today }

\author{{\bf Wen Teng }\\
{\small  School of Mathematics and Statistics, Guizhou University of Finance and Economics} \\
{\small  Guiyang  550025, P. R. of China}\\
{\small  Email: tengwen@mail.gufe.edu.cn }\\}
 \maketitle
\begin{center}
\begin{minipage}{13.cm}

{\bf \begin{center} ABSTRACT \end{center}}
In this paper, we first introduce the concept and representations of   modified $\lambda$-differential  Lie-Yamaguti algebras. We then establish the cohomology of a modified $\lambda$-differential  Lie-Yamaguti algebra with coefficients in a  representation.  As   applications, we  investigate the formal deformations and  abelian extensions of   modified $\lambda$-differential  Lie-Yamaguti algebras by using the second cohomology group.
 \smallskip

 \smallskip

 \end{minipage}
 \end{center}
 \normalsize\vskip0.5cm

\section{Introduction}
\def\theequation{\arabic{section}. \arabic{equation}}
\setcounter{equation} {0}

The concept of Lie-Yamaguti algebras  (also called a general Lie triple system or a  Lie triple algebra) was introduced by Kinyon and Weinstein \cite{K01} in the study of Courant  algebroids,
 and can be constructed from Leibniz algebras, which originated  from  Nomizu's work on the invariant affine connections on homogeneous spaces in the 1950s \cite{N54}.
Its roots can be traced back to Yamaguti's study of the general Lie triple system \cite{Y57}.
Yamaguti has been established the representation and cohomology theory of   Lie-Yamaguti algebras  in  \cite{Y67}.  In \cite{B05,B09},  Benito and his collaborators have been studied Lie-Yamaguti algebras related to simple Lie algebras and irreducible Lie-Yamaguti algebras. The  deformation and extension theory for Lie-Yamaguti algebras were investigated in \cite {L15,Z15}.  Further research on Lie-Yamaguti algebras could be found in \cite{L23,G24,S21,S22,Z22,Teng23-2,TG24} and references cited therein. All these results provide a good starting point for our further research.

The term modified $r$-matrix  stemmed from the the concept of modified classical Yang-Baxter equation, which was introduced in the work of Semenov-Tian-Shansky  \cite{Semonov}.  It was later applied
to the studies of nonabelian generalized Lax pairs and affine geometry of Lie groups. Furthermore,  Jiang and Sheng \cite{J22}  have developed cohomology and deformation theory of modified $r$-matrices.
 Motivated by the   modified $r$-matrices, in \cite{Peng},  Peng et al. introduced the concept of  modified $\lambda$-differential Lie algebras. Subsequently, the modified  algebraic structure has been widely studied in \cite{D22,L22,M22,Teng23,Teng23-3,Teng24,TG24}.

It is well known that Lie-Yamaguti algebras are a generalization of Lie algebras and Lie triple systems, it is very natural to investgate modified $\lambda$-differential  Lie-Yamaguti algebras.
This was our motivation for writing the present paper.
In precisely,  we introduce the concept of a modified $\lambda$-differential Lie Yamaguti algebra, which includes a Lie-Yamaguti algebra and a modified $\lambda$-differential operator.
 Next, we propose a representation of a  modified $\lambda$-differential  Lie-Yamaguti algebra.
In addition, we also discuss the relationship of modified $\lambda$-differential Lie-Yamaguti algebras and Lie-Yamaguti algebras with derivations, and  get some interesting results.
We define a cochain map $\Phi$, and then  the cohomology of modified $\lambda$-differential Lie Yamaguti algebras with coefficients in a representation is constructed  by using a Yamaguti coboundary operator and a cochain map $\Phi$.
 Finally,  we study  formal  deformations and  abelian extensions of a  modified $\lambda$-differential  Lie-Yamaguti algebra.
 We show that the  modified $\lambda$-differential  Lie-Yamaguti algebra is rigid if the second cohomology group is trivial.
  We prove that any abelian extension of a modified $\lambda$-differential Lie-Yamaguti algebra  has a representation and a  2-cocycle and  further prove that they are classified by the second cohomology group.

  This paper is organized as follows.
In Section \ref{sec:Representations}, we  propose  the notion of     modified $\lambda$-differential Lie-Yamaguti algebras and we introduce the representation of modified $\lambda$-differential Lie-Yamaguti algebras.
In Section \ref{sec:Cohomology},  we construct the cohomology of the modified $\lambda$-differential Lie-Yamaguti algebra.
In Section \ref{sec:deformations},  we use the cohomological approach to study formal deformations of modified $\lambda$-differential Lie-Yamaguti algebras.
 Finally, abelian extensions of modified $\lambda$-differential Lie-Yamaguti algebras are discussed in Section  \ref{sec:extensions}.

Throughout this paper, $\mathbb{K}$ denotes a field of characteristic zero. All the  vector spaces  and
   (multi)linear maps are taken over $\mathbb{K}$.

\section{ Representations of modified $\lambda$-differential Lie-Yamaguti algebras}\label{sec:Representations}
\def\theequation{\arabic{section}. \arabic{equation}}
\setcounter{equation} {0}

In this section, we introduce the   concept of  modified $\lambda$-differential Lie-Yamaguti algebras.
We explore the relationship between modified $\lambda$-differential operators and derivations,
and provide some   examples. Finally, we propose  the  representation of  modified $\lambda$-differential Lie-Yamaguti algebras.
Furthermore, we discuss its relationship with the representation of Lie-Yamaguti algebras with a derivation.

First we  recall some basic  definitions of  Lie-Yamaguti  algebra from \cite{K01} and \cite{Y67}.

\begin{definition} (\cite{K01})
A Lie-Yamaguti  algebra  is a  triple  $(\mathfrak{g}, [, ], \{, , \})$ in which $\mathfrak{g}$ is a vector space together with  a binary bracket $[, ]$ and  a ternary bracket $\{, , \}$
on $\mathfrak{g}$   satisfying
\begin{eqnarray*}
&&(LY1)~~ [x, y]=- [y,  x],\\
&&(LY2) ~~\{x,y,z\}=- \{y,x,z\},\\
&&(LY3) ~~ \circlearrowleft_{x,y,z}[[x, y], z]+\circlearrowleft_{x,y,z}\{x, y, z\}=0,\\
&&(LY4) ~~\circlearrowleft_{x,y,z}\{[x, y], z, a\}=0,\\
&&(LY5) ~~\{a, b, [x, y]\}=[\{a, b, x\}, y]+[x,\{a, b, y\}],\\
&&(LY6) ~~ \{a, b, \{x, y, z\}\}=\{\{a, b, x\}, y, z\}+ \{x,  \{a, b, y\}, z\}+ \{x,  y, \{a, b, z\}\},
\end{eqnarray*}
for any $ x, y, z, a, b\in \mathfrak{g}$  and where $\circlearrowleft_{x,y,z}$ denotes the sum over cyclic
permutation of $x,y,z$, that is $\circlearrowleft_{x,y,z}[[x, y], z]=[[x, y], z]+[[z, x], y]+[[y, z], x].$
\end{definition}

\begin{example} \label{exam:Lie algebra}
 Let $(\mathfrak{g}, [, ])$  be a Lie algebra.  Define a  ternary bracket $\{, , \}$ on $\mathfrak{g}$ by
 $$\{x,y,z\}=[[x,y],z], \forall x,y,z\in \mathfrak{g}.$$
 Then, $(\mathfrak{g}, [, ], \{, , \})$ is a   Lie-Yamaguti  algebra.
\end{example}

\begin{example} \cite{K01}
 Let $(\mathfrak{g}, [, ])$  be a Lie algebra  with a reductive decomposition $\mathfrak{g}=\mathfrak{n}\oplus\mathfrak{m}$, i.e.  $[\mathfrak{n},\mathfrak{n}]\subseteq \mathfrak{n}$ and $[\mathfrak{n},\mathfrak{m}]\subseteq \mathfrak{m}$.
Define a bilinear bracket $[,]_\mathfrak{m}$ and a trilinear bracket $\{,,\}_\mathfrak{m}$ on $\mathfrak{m}$ by
the projections of the Lie bracket:
$$[x,y]_\mathfrak{m}=\pi_\mathfrak{m}([x,y]),\{x,y,z\}_\mathfrak{m}=[\pi_\mathfrak{n}([x,y]),z], \forall x,y,z\in \mathfrak{m},$$
where $\pi_\mathfrak{n}:\mathfrak{g}\rightarrow \mathfrak{n},\pi_\mathfrak{m}:\mathfrak{g}\rightarrow \mathfrak{m}$ are the projection map.  Then, $(\mathfrak{m}, [-,-]_\mathfrak{m}, \{-, -, -\}_\mathfrak{m})$ is a   Lie-Yamaguti  algebra.
\end{example}

\begin{example}  \label{exam:Leibniz algebra} \cite{K01}
 Let $(\mathfrak{g}, \star)$  be a Leibniz algebra.   Define a binary and ternary bracket on $\mathfrak{g}$ by
 $$[x,y]=x\star y-y\star x, \{x,y,z\}=-(x\star y)\star z, \forall x,y,z\in \mathfrak{g}.$$
 Then, $(\mathfrak{g}, [, ], \{, , \})$ is a   Lie-Yamaguti  algebra.
\end{example}

\begin{example} \label{exam:2-dimensional LY algebra}
Let $\mathfrak{g}$ be a 2-dimensional  vector space  with a basis $\varepsilon_1,\varepsilon_2$.  If we define  a binary non-zero bracket $[,]$ and a ternary non-zero bracket $\{, , \}$ on $\mathfrak{g}$ as follows:
$$[\varepsilon_1,\varepsilon_2]=-[\varepsilon_2,\varepsilon_1]=\varepsilon_1, \{\varepsilon_1,\varepsilon_2, \varepsilon_2\}=-\{\varepsilon_2,\varepsilon_1, \varepsilon_2\}=\varepsilon_1,$$
then $(\mathfrak{g}, [, ], \{, , \})$ is a   Lie-Yamaguti  algebra.
\end{example}

\begin{example} \label{exam:3-dimensional LY algebra}
Let $\mathfrak{g}$ be a 3-dimensional  vector space  with a basis $\varepsilon_1$, $\varepsilon_2, \varepsilon_3$.  If we define  a binary non-zero bracket $[,]$ and a ternary non-zero bracket $\{, , \}$ on $\mathfrak{g}$ as follows:
$$[\varepsilon_1,\varepsilon_2]=-[\varepsilon_2,\varepsilon_1]=\varepsilon_3, \{\varepsilon_1,\varepsilon_2, \varepsilon_1\}=-\{\varepsilon_2,\varepsilon_1, \varepsilon_1\}=\varepsilon_3,$$
then $(\mathfrak{g}, [, ], \{, , \})$ is a   Lie-Yamaguti  algebra.
\end{example}

Next we introduce the   concept of  modified $\lambda$-differential Lie-Yamaguti algebras motivated by the   modified $\lambda$-differential Lie algebras \cite{Peng}

\begin{definition}
Let  $(\mathfrak{g}, [,], \{,,\})$ be a Lie-Yamaguti  algebra, for $\lambda\in \mathbb{K}$, a linear map $\varphi:\mathfrak{g}\rightarrow \mathfrak{g}$ is said to  a  modified $\lambda$-differential operator  if $\varphi$ satisfies
\begin{align}
 \varphi[x, y]=&[\varphi (x),y]+[x,\varphi (y)]+\lambda[x,y]\label{3.1}\\
\varphi\{x,y,z\}=&\{\varphi (x),y,z\}+\{x,\varphi (y),z\}+\{x,y,\varphi (z)\}+2\lambda\{x,y,z\},\label{3.2}
\end{align}
for any $ x, y, z\in \mathfrak{g}$.
\end{definition}

\begin{definition}
 A modified $\lambda$-differential Lie-Yamaguti algebra is a quadruple $(\mathfrak{g}, [, ], \{, , \},\varphi)$ consisting of a Lie-Yamaguti  algebra
$(\mathfrak{g}, [, ], \{, , \})$ together with a modified $\lambda$-differential operator $\varphi$.
\end{definition}

\begin{definition}
A homomorphism between two  modified $\lambda$-differential Lie-Yamaguti algebras $(\mathfrak{g}, [, ], $ $\{, , \},\varphi)$ and $(\mathfrak{g}', [, ]',  \{, , \}',\varphi')$ is a linear map $\eta: \mathfrak{g}\rightarrow \mathfrak{g}'$ satisfying
\begin{eqnarray*}
&&\eta([x, y])=[\eta(x), \eta(y)]', ~\eta(\{x, y, z\})=\{\eta(x), \eta(y),\eta(z)\}' \   \mathrm{and}\    \varphi'(\eta(x))=\eta( \varphi(x)),
\end{eqnarray*}
 for all $x,y,z\in \mathfrak{g}$. Furthermore,   $\eta$ is called an isomorphism from $\mathfrak{g}$ to $\mathfrak{g}'$ if the map $\eta$ is non-degenerate.
\end{definition}

\begin{remark}
(i) Let $\varphi$ be a modified  $\lambda$-differential operator on $(\mathfrak{g}, [, ], \{, , \})$. If $\lambda=0$, then $\varphi$ is a derivation
on $\mathfrak{g}$, and then $(\mathfrak{g}, [, ], \{, , \},\varphi)$ becomes a Lie-Yamaguti algebra with a derivation. See \cite{Guo23} for more details about  Lie-Yamaguti algebras with a derivation.

(ii) When a   Lie-Yamaguti algebra reduces  to a   Lie triple system, that is $[, ]=0$,  we get the notion of a modified
$\lambda$-differential  Lie triple system immediately.  See \cite{Teng23} for more details about  modified  $\lambda$-differential Lie  triple systems.

(iii)  When a Lie-Yamaguti algebra reduces to a Lie algebra, that is  $\{, , \}=0$, we get the notion of a  modified
$\lambda$-differential  Lie  algebra. See \cite{Peng,Teng24} for more details about  modified $\lambda$-differential  Lie  algebras.
\end{remark}

\begin{example}
Let $(\mathfrak{g}, [, ], \{, , \})$ be the 2-dimensional Lie-Yamaguti algebra given in $\mathrm{Example}$ \ref{exam:2-dimensional LY algebra}.
Then, for $k,k_1\in \mathbb{K},$
$$\varphi=\left(
        \begin{array}{cc}
          k & k_1  \\
          0 & -\lambda
        \end{array}
      \right)$$
is a modified $\lambda$-differential operator on $(\mathfrak{g}, [, ], \{, , \})$.
\end{example}

\begin{example}
Let $(\mathfrak{g}, [, ], \{, , \})$ be the 3-dimensional Lie-Yamaguti algebra given in $\mathrm{Example}$ \ref{exam:3-dimensional LY algebra}.
Then, $k,k_1,k_2,k_3,k_4\in \mathbb{K},$
$$\varphi=\left(
        \begin{array}{ccc}
          -\lambda & k_1 & 0  \\
          k_2 & k  & 0 \\
          k_3 & k_4  & k
        \end{array}
      \right)$$
is a modified  $\lambda$-differential operator on $(\mathfrak{g}, [, ], \{, , \})$.
\end{example}

\begin{example}
An identity map $\mathrm{id}_\mathfrak{g}:\mathfrak{g}\rightarrow \mathfrak{g}$   is a modified  $(-1)$-differential operator.
\end{example}

\begin{proposition}
Let  $(\mathfrak{g}, [, ], \{, , \})$ be a Lie-Yamaguti  algebra. A linear map $\varphi$ is    a  modified $\lambda$-differential operator
 if and only if  the map $\varphi+\lambda\mathrm{id}_\mathfrak{g}$ is    a  derivation on  $(\mathfrak{g}, [, ], \{, , \})$.
\end{proposition}

\begin{proof}
Equations \eqref{3.1} and \eqref{3.2}    are equivalent to
\begin{align*}
 (\varphi+\lambda\mathrm{id}_\mathfrak{g})[x, y]=&[(\varphi+\lambda\mathrm{id}_\mathfrak{g}) (x),y]+[x,(\varphi+\lambda\mathrm{id}_\mathfrak{g}) (y)],\\
 (\varphi+\lambda\mathrm{id}_\mathfrak{g})\{x,y,z\}=&\{(\varphi+\lambda\mathrm{id}_\mathfrak{g}) (x),y,z\}+\{x,(\varphi+\lambda\mathrm{id}_\mathfrak{g})(y),z\}+\{x,y,(\varphi+\lambda\mathrm{id}_\mathfrak{g}) (z)\}.
\end{align*}
The proposition follows.
\end{proof}

\begin{definition} (\cite{Y67})
Let $(\mathfrak{g}, [, ], \{, , \})$  be a  Lie-Yamaguti   algebra   and $V$ be a vector space. A representation of $(\mathfrak{g}, [, ], \{, , \})$ on $V$ consists of a linear map  $\rho: \mathfrak{g}\rightarrow \mathrm{End}(V)$ and two bilinear maps $D,\theta:\mathfrak{g}\times \mathfrak{g}\rightarrow \mathrm{End}(V)$ such that
\begin{eqnarray*}
&&(R1)~D(x,y)-\theta(y,x)+\theta(x,y)+\rho([x,y])-\rho(x)\rho(y)+\rho(y)\rho(x)=0,\\
&&(R2)~D([x,y],z)+D([y,z],x)+D([z,x],y)=0,\\
&&(R3)~\theta([x,y],a) =\theta(x,a)\rho(y)-\theta(y,a)\rho(x),\\
&&(R4)~D(a,b)\rho(x)=\rho(x)D(a,b)+\rho(\{a,b,x\}),\\
&&(R5)~\theta(x,[a,b])=\rho(a)\theta(x,b)-\rho(b)\theta(x,a),\\
&&(R6)~D(a,b)\theta(x,y)=\theta(x,y)D(a,b)+\theta(\{a,b,x\},y)+\theta(x,\{a,b,y\}),\\
&&(R7)~\theta(a,\{x,y,z\}) =\theta(y,z)\theta(a,x)-\theta(x,z)\theta(a,y)+D(x,y)\theta(a,z),
\end{eqnarray*}
for all  $x,y,z,a,b\in \mathfrak{g}$. In this case, we also call $V$ a $\mathfrak{g}$-module. We denote a representation by $(V; \rho, \theta, D)$.
\end{definition}
It can be concluded from (R6) that
\begin{eqnarray*}
&&(R6)'~D(a,b)D(x,y)=D(x,y)D(a,b)+D(\{a,b,x\},y)+D(x,\{a,b,y\}).
\end{eqnarray*}

\begin{example}
Let $(\mathfrak{g}, [, ], \{, , \})$  be a  Lie-Yamaguti   algebra.  We   define linear maps $\mathrm{ad}: \mathfrak{g}\rightarrow \mathrm{End}(\mathfrak{g}),\mathcal{L},\mathcal{R}:\otimes^2 \mathfrak{g}\rightarrow \mathrm{End}(\mathfrak{g})$ by
\begin{eqnarray*}
& \mathrm{ad}(x)(z):=[x,z],\mathcal{L}(x,y)(z):=\{x,y,z\},\mathcal{R}(x,y)(z):=\{z,x,y\},
\end{eqnarray*}for all $x,y,z\in \mathfrak{g}$.  Then $(\mathfrak{g};\mathrm{ad},\mathcal{L},\mathcal{R})$ forms a representation of $\mathfrak{g}$ on itself,  called the adjoint representation.
\end{example}

\begin{definition}
  A representation of the modified $\lambda$-differential Lie-Yamaguti algebra $(\mathfrak{g}, [, ],  \{, , \},\varphi)$  is a quintuple $(V; \rho,   \theta, D, \varphi_V)$ such that the following conditions are satisfied:

  (i) $(V; \rho,   \theta, D)$ is a representation of the  Lie-Yamaguti algebra $(\mathfrak{g}, [, ],  \{, , \})$;

  (ii)  $\varphi_V:V\rightarrow V$  is a linear map satisfying the following equations
\begin{align}
\varphi_V(\rho(x)u)=&\rho(\varphi(x))u+\rho(x)\varphi_V(u)+\lambda\rho(x)u,\label{3.3}\\
\varphi_V(\theta(x, y)u)=&\theta(\varphi(x), y) u +\theta(x, \varphi(y))u+\theta(x, y)\varphi_V(u)+2\lambda\theta(x, y)u,\label{3.4}
\end{align}
for any $x,y\in \mathfrak{g}$ and $ u\in V.$
\end{definition}

It can be concluded from  Eq. \eqref{3.4}  that
\begin{align}
\varphi_V(D(x, y)u)=&D(\varphi(x), y) u +D(x, \varphi(y))u+D(x, y)\varphi_V(u)+2\lambda D(x, y)u.\label{3.5}
\end{align}

\begin{example}
$(\mathfrak{g};\mathrm{ad},\mathcal{L},\mathcal{R},\varphi)$ is an adjoint representation of the modified $\lambda$-differential  Lie-Yamaguti algebra $(\mathfrak{g},  [, ],  \{, , \}, \varphi)$.
\end{example}

\begin{remark}
 Let $(V; \rho,   \theta, D, \varphi_V)$ be a representation of the modified $\lambda$-differential Lie-Yamaguti algebra $(\mathfrak{g}, $ $[, ],  \{, , \},\varphi)$. If $\lambda=0$, then $(V; \rho,   \theta, D, \varphi_V)$ be a representation of a Lie-Yamaguti algebra  with a derivation $(\mathfrak{g}, [, ],  \{, , \},\varphi)$.
\end{remark}

The representation of  modified  $\lambda$-differential Lie-Yamaguti algebras is closely related to the representation of  a Lie-Yamaguti algebra  with a derivation .

\begin{proposition}
The quintuple $(V; \rho,   \theta, D, \varphi_V)$ is a representation of the modified $\lambda$-differential Lie-Yamaguti algebra $(\mathfrak{g}, $ $[, ],  \{, , \},\varphi)$
 if and only if  $(V; \rho,   \theta, D, \varphi_V+\lambda \mathrm{id}_V)$ be a representation of a Lie-Yamaguti algebra  with a derivation $(\mathfrak{g}, [, ],  \{, , \},\varphi+\lambda \mathrm{id}_{\mathfrak{g}})$.
\end{proposition}

\begin{proof}
Equations \eqref{3.3}-\eqref{3.5}    are equivalent to
\begin{align*}
(\varphi_V+\lambda \mathrm{id}_V)(\rho(x)u)=&\rho((\varphi+\lambda \mathrm{id}_{\mathfrak{g}})(x))u+\rho(x)(\varphi_V+\lambda \mathrm{id}_V)(u),\\
(\varphi_V+\lambda \mathrm{id}_V)(\theta(x, y)u)=&\theta((\varphi+\lambda \mathrm{id}_{\mathfrak{g}})(x), y) u +\theta(x, (\varphi+\lambda \mathrm{id}_{\mathfrak{g}})(y))u+\theta(x, y)(\varphi_V+\lambda \mathrm{id}_V)(u),\\
(\varphi_V+\lambda \mathrm{id}_V)(D(x, y)u)=&D((\varphi+\lambda \mathrm{id}_{\mathfrak{g}})(x), y) u +D(x, (\varphi+\lambda \mathrm{id}_{\mathfrak{g}})(y))u+D(x, y)(\varphi_V+\lambda \mathrm{id}_V)(u).
\end{align*}
The proposition follows.
\end{proof}

Next we  construct the semidirect product of the modified $\lambda$-differential  Lie-Yamaguti algebra.

\begin{proposition}
If  $(V; \rho,   \theta, D, \varphi_V)$  is a representation of the  modified  $\lambda$-differential  Lie-Yamaguti algebra $(\mathfrak{g}, [, ],  \{, , \}, \varphi)$, then $\mathfrak{g} \oplus V$ is a modified  $\lambda$-differential  Lie-Yamaguti   algebra  under the following maps:
\begin{eqnarray*}
&&[x+u, y+v]_{\ltimes}:=[x, y]+\rho(x)(v)-  \rho(y)(u),\\
&&\{x+u, y+v, z+w\}_{\ltimes}:=\{x, y, z\}+D(x, y)(w)- \theta(x, z)(v)+ \theta(y, z)(u),\\
&&\varphi\oplus \varphi_V(x+u)=\varphi(x)+\varphi_V(u),
\end{eqnarray*}
for all  $x, y, z\in \mathfrak{g}$ and $u, v, w\in V$. In the case, the modified $\lambda$-differential  Lie-Yamaguti   algebra $\mathfrak{g} \oplus V$  is called a semidirect product of $\mathfrak{g}$ and $V$, denoted by $\mathfrak{g}\ltimes V=(\mathfrak{g} \oplus V,[,]_{\ltimes},\{,,\}_{\ltimes},\varphi\oplus \varphi_V)$.
\end{proposition}

\begin{proof}
In view of \cite{Y67}, $(\mathfrak{g} \oplus V,[,]_{\ltimes},\{,,\}_{\ltimes})$ is a Lie-Yamaguti algebra.

Next, for any $x, y, z\in \mathfrak{g}$ and $u, v, w\in V$,   in view of  Eqs. \eqref{3.1}-- \eqref{3.5},  we have
\begin{align*}
&\varphi\oplus \varphi_V[x+u, y+v]_{\ltimes}\\
&=\varphi\oplus \varphi_V([x, y]+\rho(x)(v)-  \rho(y)(u))\\
&=\varphi[x, y]+\varphi_V(\rho(x)(v)-  \rho(y)(u))\\
&=[\varphi (x),y]+[x,\varphi (y)]+\lambda[x,y]+\rho(\varphi(x))v+\rho(x)\varphi_V(v)+\lambda\rho(x)v- \rho(\varphi(y))u-\rho(y)\varphi_V(u)\\
&~~~~-\lambda\rho(y)u\\
&=\big([\varphi (x),y]+\rho(\varphi(x))v-\rho(y)\varphi_V(u)\big)+\big([x,\varphi (y)]+\rho(x)\varphi_V(v)- \rho(\varphi(y))u\big)\\
&~~~~+\lambda([x,y]+\rho(x)v-\rho(y)u)\\
&=[\varphi\oplus \varphi_V (x+u),y+v]+[x+u,\varphi\oplus \varphi_V (y+v)]+\lambda[x+u,y+v],\\
&\varphi\oplus \varphi_V\{x+u, y+v, z+w\}_{\ltimes}\\
&=\varphi\oplus \varphi_V(\{x, y, z\}+D(x, y)(w)- \theta(x, z)(v)+ \theta(y, z)(u))\\
&=\varphi\{x, y, z\}+\varphi_V(D(x, y)(w)- \theta(x, z)(v)+ \theta(y, z)(u))\\
&=\{\varphi (x),y,z\}+\{x,\varphi (y),z\}+\{x,y,\varphi (z)\}+2\lambda\{x,y,z\}+D(\varphi(x), y)w +D(x, \varphi(y))w\\
&~~~~+D(x, y)\varphi_V(w)+2\lambda D(x, y)w- \theta(\varphi(x), z) v -\theta(x, \varphi(z))v-\theta(x,z)\varphi_V(v)-2\lambda\theta(x, z)v\\
&~~~~+ \theta(\varphi(y), z) u +\theta(y, \varphi(z))u+\theta(y, z)\varphi_V(u)+2\lambda\theta(y, z)u\\
&=\big(\{\varphi (x),y,z\}+D(\varphi(x), y)w - \theta(\varphi(x), z) v+\theta(y, z)\varphi_V(u)\big)+\big(\{x,\varphi (y),z\}+D(x, \varphi(y))w\\
&~~~~-\theta(x,z)\varphi_V(v)+ \theta(\varphi(y), z) u\big)+\big(\{x,y,\varphi (z)\}+D(x, y)\varphi_V(w) -\theta(x, \varphi(z))v+\theta(y, \varphi(z))u\big)\\
&~~~~+2\lambda\big(\{x,y,z\}+ D(x, y)w-\theta(x, z)v+\theta(y, z)u\big)\\
&=\{\varphi\oplus \varphi_V(x+u), y+v, z+w\}_{\ltimes}+\{x+u, \varphi\oplus \varphi_V(y+v), z+w\}_{\ltimes}\\
&~~~~+\{x+u, y+v, \varphi\oplus \varphi_V(z+w)\}_{\ltimes}+2\lambda\{x+u, y+v, z+w\}_{\ltimes}.
\end{align*}
So  $(\mathfrak{g} \oplus V,[,]_{\ltimes},\{,,\}_{\ltimes},\varphi\oplus \varphi_V)$ is a modified  $\lambda$-differential  Lie-Yamaguti   algebra.
\end{proof}

\section{ Cohomology of modified $\lambda$-differential  Lie-Yamaguti algebras}\label{sec:Cohomology}
\def\theequation{\arabic{section}. \arabic{equation}}
\setcounter{equation} {0}

In this section, we will construct the cohomology of  modified $\lambda$-differential  Lie-Yamaguti algebras.
In the next section, we will  use the second cohomology group to study the formal deformation and abelian extension of   modified $\lambda$-differential  Lie-Yamaguti algebras.

 Firstly, let us recall the Yamaguti cohomology theory on Lie-Yamaguti algebras in \cite{Y67}.  Let $(V; \rho, \theta, D)$ be a
representation of a Lie-Yamaguti algebra $(\mathfrak{g}, [, ], \{, , \})$, and we denote the set of $(n+1)$-cochains by
$\mathcal{C}^{n+1}_{\mathrm{LY}}(\mathfrak{g},V)$, where
\begin{equation*}
\mathcal{C}_{\mathrm{LY}}^{n+1}(\mathfrak{g},V)= \left\{ \begin{array}{ll}
\mathrm{Hom}(\underbrace{\wedge^2 \mathfrak{g}\otimes\cdots\otimes\wedge^2 \mathfrak{g}}_n,V)\times \mathrm{Hom}(\underbrace{\wedge^2 \mathfrak{g}\otimes\cdots\otimes\wedge^2 \mathfrak{g}}_n\otimes \mathfrak{g},V) &\mbox{ \mbox{}  $n\geq 1,$  }\\
$$\mathrm{Hom}(\mathfrak{g},V)$$ &\mbox{ \mbox{}  $ n=0$.  }
 \end{array}  \right.
\end{equation*}

In the sequel, we recall the coboundary map of $(n+1)$-cochains on a Lie-Yamaguti algebra $\mathfrak{g}$ with the
coefficients in the representation $(V; \rho, \theta, D)$:

If $n\geq 1$, for any $(f,g)\in \mathcal{C}_{\mathrm{LY}}^{n+1}(\mathfrak{g},V)$, $ \mathfrak{K}_i=x_i\wedge y_i\in \wedge^2 \mathfrak{g}, (i=1,2,\cdots,n+1), z\in \mathfrak{g}$,  the coboundary map $\delta^{n+1}=(\delta^{n+1}_I,\delta^{n+1}_{II}):\mathcal{C}_{\mathrm{LY}}^{n+1}(\mathfrak{g},V)\rightarrow \mathcal{C}_{\mathrm{LY}}^{n+2}(\mathfrak{g},V), (f,g)\mapsto (\delta^{n+1}_I(f,g),\delta^{n+1}_{II}(f,g))$ is given as follows:
\begin{align*}
&\delta^{n+1}_I(f,g)(\mathfrak{K}_1,\cdots,\mathfrak{K}_{n+1})\\
=&(-1)^n(\rho(x_{n+1})g(\mathfrak{K}_1,\cdots,\mathfrak{K}_{n},y_{n+1})-\rho(y_{n+1})g(\mathfrak{K}_1,\cdots,\mathfrak{K}_{n},x_{n+1})\\
&-g(\mathfrak{K}_1,\cdots,\mathfrak{K}_{n},[x_{n+1},y_{n+1}]))+\sum_{k=1}^{n}(-1)^{k+1}D(\mathfrak{K}_k)f(\mathfrak{K}_1,\cdots,\widehat{\mathfrak{K}_{k}}\cdots,\mathfrak{K}_{n+1})\\
&+\sum_{1\leq k<l\leq n+1}(-1)^k f(\mathfrak{K}_1,\cdots,\widehat{\mathfrak{K}_{k}}\cdots,\{x_k,y_k,x_l\}\wedge y_l+x_l\wedge \{x_k,y_k,y_l\},\cdots,\mathfrak{K}_{n+1}),
\end{align*}
\begin{align*}
&\delta^{n+1}_{II}(f,g)(\mathfrak{K}_1,\cdots,\mathfrak{K}_{n+1},z)\\
=&(-1)^n(\theta(y_{n+1},z)g(\mathfrak{K}_1,\cdots,\mathfrak{K}_{n},x_{n+1})-\theta(x_{n+1},z)g(\mathfrak{K}_1,\cdots,\mathfrak{K}_{n},y_{n+1}))\\
&+\sum_{k=1}^{n+1}(-1)^{k+1}D(\mathfrak{K}_k)g(\mathfrak{K}_1,\cdots,\widehat{\mathfrak{K}_{k}}\cdots,\mathfrak{K}_{n+1},z)\\
&+\sum_{1\leq k<l\leq n+1}(-1)^k g(\mathfrak{K}_1,\cdots,\widehat{\mathfrak{K}_{k}}\cdots,\{x_k,y_k,x_l\}\wedge y_l+x_l\wedge \{x_k,y_k,y_l\},\cdots,\mathfrak{K}_{n+1},z)\\
&+\sum_{k=1}^{n+1}(-1)^kg(\mathfrak{K}_1,\cdots,\widehat{\mathfrak{K}_{k}}\cdots,\mathfrak{K}_{n+1},\{x_k,y_k,z\}).
\end{align*}
where $~\widehat{}~$ denotes omission.
For the case that $n=0$, for any $f\in \mathcal{C}_{\mathrm{LY}}^1(\mathfrak{g},V)$ , the coboundary map
 $\delta^1=(\delta^1_I,\delta^1_{II})\text{:}$\, $\mathcal{C}_{\mathrm{LY}}^1(\mathfrak{g},V)\rightarrow \mathcal{C}_{\mathrm{LY}}^2(\mathfrak{g},V),f\rightarrow (\delta^1_I(f),\delta^1_{II}(f))$ is given by:
\begin{align*}
\delta^1_I(f)(x,y)=&\rho(x)f(y)-\rho(y)f(x)-f([x,y]),\\
\delta^1_{II}(f)(x,y,z)=&D(x,y)f(z)+\theta(y,z)f(x)-\theta(x,z)f(y)-f(\{x,y,z\}).
\end{align*}
The corresponding cohomology groups are denoted by $\mathcal{H}^{\bullet}_{\mathrm{LY}}(\mathfrak{g},V).$

Next  we   construct a cohomology theory for modified $\lambda$-differential Lie-Yamaguti algebras by using the Yamaguti cohomology.

For any  $(f,g)\in \mathcal{C}_{\mathrm{LY}}^{n+1}(\mathfrak{g},V)$, $n\geq 1$, we define a linear map $\Phi^{n+1}=(\Phi_I^{n+1},\Phi^{n+1}_{II}):$
$\mathcal{C}^{n+1}_{\mathrm{LY}}(\mathfrak{g},V)\rightarrow \mathcal{C}^{n+1}_{\mathrm{LY}}(\mathfrak{g},V),  (f,g)\mapsto (\Phi_I^{n+1}(f), \Phi_{II}^{n+1}(g))$  by:
\begin{align*}
\Phi^{n+1}_{I}(f)=&\sum_{i=1}^{2n}f\circ(\mathrm{Id}^{i-1},\varphi,\mathrm{Id}^{2n-i})\circ((\mathrm{Id}\wedge \mathrm{Id})^{n})+(2n-1) f\circ((\mathrm{Id}\wedge \mathrm{Id})^{n}) -\varphi_V\circ f\circ((\mathrm{Id}\wedge \mathrm{Id})^{n}),\\
\Phi^{n+1}_{II}(g)=&\sum_{i=1}^{2n+1}g\circ(\mathrm{Id}^{i-1},\varphi,\mathrm{Id}^{2n+1-i})\circ((\mathrm{Id}\wedge \mathrm{Id})^{n}\wedge \mathrm{Id})+2n g\circ((\mathrm{Id}\wedge \mathrm{Id})^{n}\wedge \mathrm{Id}) \\
&\ ~~~-\varphi_V\circ g\circ((\mathrm{Id}\wedge \mathrm{Id})^{n}\wedge \mathrm{Id}),
\end{align*}
In particular, when $n=0$, define $\Phi^{1}:\mathcal{C}^{1}_{\mathrm{LY}}(\mathfrak{g},V)\rightarrow \mathcal{C}^{1}_{\mathrm{LY}}(\mathfrak{g},V)$ by
$$\Phi^1(f)=f\circ \varphi-\varphi_V\circ f.$$

\begin{lemma}\label{lemma:chain map}
The map $\Phi^{n+1}:\mathcal{C}^{n+1}_{\mathrm{LY}}(\mathfrak{g},V)\rightarrow \mathcal{C}^{n+1}_{\mathrm{LY}}(\mathfrak{g},V)$ is a cochain map, that is, $\delta^{n}\circ\Phi^{n+1}=\Phi^{n+1}\circ\delta^n$
\end{lemma}
\begin{proof}
It follows by a straightforward tedious calculations.
\end{proof}

\begin{definition}
Let $(V; \rho,   \theta, D, \varphi_V)$  be a representation of the  modified  $\lambda$-differential  Lie-Yamaguti algebra $(\mathfrak{g},  [, ],   \{, , \}, \varphi)$.
We define the set of modified  $\lambda$-differential  Lie-Yamaguti algebra 1-cochains to be $\mathcal{C}^{1}_{\mathrm{MDLY}}(\mathfrak{g},V)=\mathcal{C}^{1}_{\mathrm{LY}}(\mathfrak{g},V)$.
  For $n\geq 1$, we define the set of modified  $\lambda$-differential  Lie-Yamaguti algebra $(n+1)$-cochains by
$$\mathcal{C}_{\mathrm{MDLY}}^{n+1}(\mathfrak{g},V):=
\mathcal{C}^{n+1}_{\mathrm{LY}}(\mathfrak{g},V)\oplus \mathcal{C}^{n}_{\mathrm{LY}}(\mathfrak{g},V).$$

Define  a linear map  $\partial^1:\mathcal{C}_{\mathrm{MDLY}}^{1}(\mathfrak{g},V)\rightarrow \mathcal{C}_{\mathrm{MDLY}}^{2}(\mathfrak{g},V)$  by
\begin{align*}
\partial^1(f_1)=(\delta^1 f_1, -\Phi^1 (f_1)), \forall f_1\in \mathcal{C}_{\mathrm{MDLY}}^{1}(\mathfrak{g},V).
\end{align*}
For  $n = 1$,   we define  the linear map  $\partial^{2}:\mathcal{C}_{\mathrm{MDLY}}^{2}(\mathfrak{g},V)\rightarrow \mathcal{C}_{\mathrm{MDLY}}^{3}(\mathfrak{g},V)$  by
\begin{align*}
\partial^{2}((f_1,g_1),f_2)=(\delta^{2}(f_1,g_1), \delta^{1} (f_2)+ \Phi^{2}(f_1,g_1)),~ \forall~ ((f_1,g_1),f_2)\in \mathcal{C}_{\mathrm{MDLY}}^{2}(\mathfrak{g},V)
\end{align*}
Then, for  $n \geq 2$,   we define  the linear map  $\partial^{n+1}:\mathcal{C}_{\mathrm{MDLY}}^{n+1}(\mathfrak{g},V)\rightarrow \mathcal{C}_{\mathrm{MDLY}}^{n+2}(\mathfrak{g},V)$  by
\begin{align*}
\partial^{n+1}((f_1,g_1),(f_2,g_2))=(\delta^{n+1}(f_1,g_1), \delta^{n} (f_2,g_2)+(-1)^{n+1} \Phi^{n+1}(f_1,g_1)),
\end{align*}
for any  $((f_1,g_1),(f_2,g_2))\in \mathcal{C}_{\mathrm{MDLY}}^{n+1}(\mathfrak{g},V).$
\end{definition}

In view of Lemma \ref{lemma:chain map}, we have the following theorem.

\begin{theorem}\label{theorem: cochain complex}
The map $\partial^{n+1}$ is a coboundary operator, i.e., $\partial^{n+1}\circ\partial^{n}=0.$
\end{theorem}

Therefore, from Theorem \ref{theorem: cochain complex}, we get a cochain complex $(\mathcal{C}_{\mathrm{MDLY}}^{\bullet}(\mathfrak{g},V),\partial^{\bullet})$   called the cochain complex
of modified $\lambda$-differential  Lie-Yamaguti algebra $(\mathfrak{g},  [, ],   \{, , \}, \varphi)$ with coefficients in $(V; \rho,   \theta, D, \varphi_V)$.
The cohomology of  $(\mathcal{C}_{\mathrm{MDLY}}^{\bullet}(\mathfrak{g},V),\partial^{\bullet})$, denoted by
 $\mathcal{H}_{\mathrm{MDLY}}^{\bullet}(\mathfrak{g},V)$,
 is called the cohomology of the  modified  $\lambda$-differential   Lie-Yamaguti algebra  $(\mathfrak{g},  [, ],   \{, , \}, \varphi)$ with coefficients in $(V; \rho,   \theta, D, \varphi_V)$.

 In particular,  when $(V; \rho,   \theta, D, \varphi_V)=(\mathfrak{g};\mathrm{ad},\mathcal{L},\mathcal{R},\varphi)$, we just denote $(\mathcal{C}_{\mathrm{MDLY}}^{\ast}(\mathfrak{g},\mathfrak{g}),\partial^{\bullet})$, $\mathcal{H}_{\mathrm{MDLY}}^{\bullet}(\mathfrak{g},\mathfrak{g})$  by
 $(\mathcal{C}_{\mathrm{MDLY}}^{\bullet}(\mathfrak{g}),\partial^{\bullet})$, $\mathcal{H}_{\mathrm{MDLY}}^{\bullet}(\mathfrak{g})$ respectively, and call them the cochain complex, the cohomology of modified  $\lambda$-differential   Lie-Yamaguti algebra $(\mathfrak{g},  [, ], \{, , \}, \varphi)$ respectively.

It is obvious that there is a  short exact sequence of cochain complexes:
\begin{align*}
0\rightarrow \mathcal{C}_{\mathrm{LY}}^{\bullet-1}(\mathfrak{g},V)\stackrel{}{\longrightarrow}\mathcal{C}_{\mathrm{MDLY}}^{\bullet}(\mathfrak{g},V)\stackrel{}{\longrightarrow}\mathcal{C}_{\mathrm{LY}}^{\bullet}(\mathfrak{g},V)\rightarrow 0.
\end{align*}
It induces a long exact sequence of cohomology groups:
\begin{align*}
\cdots\rightarrow \mathcal{H}_{\mathrm{MDLY}}^{p}(\mathfrak{g},V)\rightarrow \mathcal{H}_{\mathrm{LY}}^{p}(\mathfrak{g},V)\rightarrow \mathcal{H}_{\mathrm{MDLY}}^{p+1}(\mathfrak{g},V)\rightarrow \mathcal{H}_{\mathrm{LY}}^{p+1}(\mathfrak{g},V)\rightarrow \cdots.
\end{align*}

 \section{  Formal deformations of modified $\lambda$-differential Lie-Yamaguti algebras}\label{sec:deformations}
\def\theequation{\arabic{section}. \arabic{equation}}
\setcounter{equation} {0}

In this section, we study   formal deformations of modified $\lambda$-differential Lie-Yamaguti algebras.
Let $\mathbb{K}[[t]]$ be a ring of power series of one variable $t$, and let $\mathfrak{g}[[t]]$ be the set of formal power series over $\mathfrak{g}$.  If
$(\mathfrak{g}, [, ], \{, , \})$ is a Lie-Yamaguti algebra, then there is a Lie-Yamaguti algebra structure over the
ring $\mathbb{K}[[t]]$ on $\mathfrak{g}[[t]]$ given by
\begin{align*}
&[\sum_{i=0}^{\infty}x_it^i,\sum_{i=0}^{\infty}y_jt^j]=\sum_{s=0}^{\infty}\sum_{i+j=s}[x_i,y_j]t^s,\\
& \{\sum_{i=0}^{\infty}x_it^i,\sum_{i=0}^{\infty}y_jt^j,\sum_{k=0}^{\infty}z_kt^k\}=\sum_{s=0}^{\infty}\sum_{i+j+k=s}\{x_i,y_j,z_k\}t^s.
\end{align*}

\begin{definition}
A   formal deformation of the modified $\lambda$-differential  Lie-Yamaguti algebra  $(\mathfrak{g}, [, ], \{, , \},\varphi)$ is a triple  $ (f_t, g_t, \varphi_t)$  of the forms
$$f_t=\sum_{i=0}^{\infty}f_it^i,~~ g_t=\sum_{i=1}^{\infty}g_it^i,~~\varphi_t=\sum_{i=0}^{\infty}\varphi_it^i,$$
such that the following conditions are satisfied:

(i) $((f_i,g_i),\varphi_i)\in \mathcal{C}^{2}_{\mathrm{MDLY}}(\mathfrak{g});$

(ii) $f_0=[, ],g_0=\{, , \}$ and $\varphi_0=\varphi;$

(iii) and $(\mathfrak{g}[[t]], f_t, g_t, \varphi_t)$  is  a   modified $\lambda$-differential  Lie-Yamaguti algebra over $\mathbb{K}[[t]]$.
\end{definition}

Let $ (f_t, g_t, \varphi_t)$  be a formal deformation as above. Then, for any $ x, y, z, a, b\in \mathfrak{g}$, the following equations must hold:
\begin{align*}
&~f_t(x, y)+f_t(y,  x)=0,~~g_t(x,y,z)+g_t(y,x,z)=0,\\
&~ \circlearrowleft_{x,y,z}f_t(f_t(x, y), z)+\circlearrowleft_{x,y,z}g_t(x,y,z)=0,\\
&~\circlearrowleft_{x,y,z}g_t(f_t(x, y), z, a)=0,\\
&~g_t(a, b, f_t(x, y))=f_t(g_t(a, b, x), y)+f_t(x,g_t(a, b, y)),\\
& ~ g_t(a, b, g_t(x, y, z))=g_t(g_t(a, b, x), y, z)+ g_t(x,  g_t(a, b, y), z)+ g_t(x,  y, g_t(a, b, z)),\\
&~ \varphi_t(f_t(x, y))=f_t(\varphi_t(x),y)+f_t(x,\varphi_t(y))+\lambda f_t(x,y), \\
&~\varphi_t(g_t(x,y,z))=g_t(\varphi_t(x),y,z)+g_t(x,\varphi_t(y),z)+g_t(x,y,\varphi_t(z))+2\lambda g_t(x,y,z)).
\end{align*}
Collecting the coefficients of $t^n$, we get that the above equations are equivalent to the following equations.
\begin{small}
\begin{align}
&~f_n(x, y)+f_n(y,  x)=0,~~g_n(x,y,z)+g_n(y,x,z)=0,\label{4.1}\\
& \sum_{i=0}^n\circlearrowleft_{x,y,z}f_i(f_{n-i}(x, y), z)+\circlearrowleft_{x,y,z}g_n(x,y,z)=0,\label{4.2}\\
&\sum_{i=0}^n\circlearrowleft_{x,y,z}g_i(f_{n-i}(x, y), z, a)=0,\label{4.3}\\
&\sum_{i=0}^ng_i(a, b, f_{n-i}(x, y))=\sum_{i=0}^nf_i(g_{n-i}(a, b, x), y)+\sum_{i=0}^nf_i(x,g_{n-i}(a, b, y)),\label{4.4}\\
& \sum_{i=0}^n g_i(a, b, g_{n-i}(x, y, z))=\sum_{i=0}^n\big(g_i(g_{n-i}(a, b, x), y, z)+g_i(x,  g_{n-i}(a, b, y), z)+ g_i(x,  y, g_{n-i}(a, b, z))\big),\label{4.5}\\
&~ \sum_{i=0}^n\varphi_i(f_{n-i}(x, y))=\sum_{i=0}^n\big(f_i(\varphi_{n-i}(x),y)+f_i(x,\varphi_{n-i}(y))\big)+\lambda f_n(x,y), \label{4.6}\\
&\sum_{i=0}^n\varphi_i(g_{n-i}(x,y,z))=\sum_{i=0}^n\big(g_i(\varphi_{n-i}(x),y,z)+g_i(x,\varphi_{n-i}(y),z)+g_i(x,y,\varphi_{n-i}(z))\big)+2\lambda g_n(x,y,z).\label{4.7}
\end{align}
\end{small}
Note that for $n=0$, equations \eqref{4.1}--\eqref{4.7}  are equivalent to $(\mathfrak{g}, f_0, g_0, \varphi_0)$  is  a   modified  $\lambda$-differential  Lie-Yamaguti algebra.

\begin{proposition}\label{prop:2-cocycle}
Let $ (f_t, g_t, \varphi_t)$  be a formal deformation of  a modified $\lambda$-differential  Lie-Yamaguti algebra  $(\mathfrak{g}, [, ], \{, , \},\varphi)$.
Then $((f_1,g_1),\varphi_1)$ is a 2-cocycle in the cochain complex $(\mathcal{C}^{\bullet}_{\mathrm{MDLY}}(\mathfrak{g}),\partial^{\bullet})$.
\end{proposition}
\begin{proof}
For  $n=1$,   equations \eqref{4.1}--\eqref{4.5} become
\begin{small}
\begin{align}
&f_1(x, y)+f_1(y,  x)=0,~~g_1(x,y,z)+g_1(y,x,z)=0,\label{4.8}\\
&[f_{1}(x, y), z]+f_1([x, y], z)+[f_{1}(z, x), y]+f_1([z, x], y)+[f_{1}(y, z), x]+f_1([y, z], x)\nonumber\\
&+g_1(x,y,z)+g_1(z,x,y)+g_1(y,z, x)=0,\label{4.9}\\
&g_1([x, y], z, a)+\{f_{1}(x, y), z, a\}+g_1([z, x], y, a)+\{f_{1}(z, x), y, a\}+g_1([y, z], x, a)\nonumber\\
&+\{f_{1}(y, z), x, a\}=0,\label{4.10}\\
&g_1(a, b, [x, y])+\{a, b, f_{1}(x, y)\}\nonumber\\
&=f_1(\{a, b, x\}, y)+[g_{1}(a, b, x), y]+f_1(x,\{a, b, y\})+[x,g_{1}(a, b, y)],\label{4.11}\\
&  g_1(a, b, \{x, y, z\})+  \{a, b, g_{1}(x, y, z)\}=g_1(\{a, b, x\}, y, z)+ \{g_{1}(a, b, x), y, z\}\nonumber\\
&+ g_1(x,  \{a, b, y\}, z)+ \{x,  g_{1}(a, b, y), z\}+ g_1(x,  y, \{a, b, z\})+ \{x,  y, g_{1}(a, b, z)\}.\label{4.12}
\end{align}
\end{small}
From Eqs. \eqref{4.8}--\eqref{4.12}, we get $(f_1,g_1)\in \mathcal{C}^{2}_{\mathrm{LY}}(\mathfrak{g})$ and $\delta^{2}(f_1,g_1)=0$. 
For  $n=1$,   equations \eqref{4.6} and \eqref{4.7} become
\begin{small}
\begin{align}
&\varphi(f_{1}(x, y))+\varphi_1([x, y])=[\varphi_{1}(x),y]+f_1(\varphi(x),y)+[x,\varphi_{1}(y)]+f_1(x,\varphi(y)]+\lambda f_1(x,y), \label{4.13}\\
&\varphi(g_{1}(x,y,z))+\varphi_1(\{x,y,z\})=\{\varphi_{1}(x),y,z\}+g_1(\varphi(x),y,z)+\{x,\varphi_{1}(y),z\}+g_1(x,\varphi(y),z)\nonumber\\
&~~~~~~~~~~~~~~~~~~~~~~~~~~~~~~~~~~~~~~~~+\{x,y,\varphi_{1}(z)\}+g_1(x,y,\varphi(z))+2\lambda g_1(x,y,z).\label{4.14}
\end{align}
\end{small}
 Further by  Eqs. \eqref{4.13} and \eqref{4.14}, we have
$\delta_I^{1} (\varphi_1)+\Phi_I^{2}(f_1)=0$ and $\delta_{II}^{1} (\varphi_1)+ \Phi_{II}^{2}(g_1)=0$  respectively.
 Hence, $\partial^{2}((f_1,g_1),\varphi_1)=0,$ that is, $((f_1,g_1),\varphi_1)$ is a 2-cocycle in   $(\mathcal{C}^{\bullet}_{\mathrm{MDLY}}(\mathfrak{g}),\partial^\bullet)$.
\end{proof}

\begin{definition}
The 2-cocycle  $((f_1,g_1),\varphi_1)$ is called the infinitesimal of the   formal deformation $ (f_t, g_t, \varphi_t)$ of a modified $\lambda$-differential  Lie-Yamaguti algebra  $(\mathfrak{g}, [, ], \{, , \},\varphi)$.
\end{definition}

\begin{definition}
Let $(f_t, g_t, \varphi_t)$ and $(f'_t, g'_t, \varphi'_t)$ be two formal deformations of
a  modified $\lambda$-differential  Lie-Yamaguti algebra  $(\mathfrak{g}, [, ], \{, , \},\varphi)$. A formal isomorphism from $(\mathfrak{g}[[t]], f_t, g_t, \varphi_t)$ to $(\mathfrak{g}[[t]],f'_t, g'_t, \varphi'_t)$ is a power series
$\Psi_t=\sum_{i=o}^{\infty}\Psi_it^i:\mathfrak{g}[[t]]\rightarrow \mathfrak{g}[[t]]$ , where  $\Psi_i:\mathfrak{g}\rightarrow \mathfrak{g}$ are linear maps
with $\Psi_0=\mathrm{Id}_\mathfrak{g}$, such that:
\begin{align}
& \Psi_t \circ f_t=f'_t\circ(\Psi_t\otimes\Psi_t),\label{4.15}\\
& \Psi_t \circ g_t=g'_t\circ(\Psi_t\otimes\Psi_t\otimes\Psi_t),\label{4.16}\\
& \Psi_t \circ \varphi_t=\varphi'_t\circ \Psi_t,\label{4.17}
\end{align}
In this case, we say that the two   formal deformations $(f_t, g_t, \varphi_t)$ and $(f'_t, g'_t, \varphi'_t)$
are equivalent.
\end{definition}

\begin{proposition}
The infinitesimals of two equivalent   formal deformations of $(\mathfrak{g}, [, ], \{, , \},\varphi)$
are in the same cohomology class in $\mathcal{H}^{2}_{\mathrm{MDLY}}(\mathfrak{g})$.
\end{proposition}
\begin{proof}
Let $\Psi_t:(\mathfrak{g}[[t]], f_t, g_t, \varphi_t)\rightarrow (\mathfrak{g}[[t]],f'_t, g'_t, \varphi'_t)$ be a formal isomorphism.
 By expanding Eqs. \eqref{4.15}--\eqref{4.17} and comparing the coefficients of $t$ on both sides, we have
 \begin{align*}
&f_1-f'_1=f_0\circ (\Psi_1\o \mathrm{Id}_\mathfrak{g})+f_0\circ (\mathrm{Id}_\mathfrak{g} \o\Psi_1)-\Psi_1\circ f_0,\\
&g_1-g'_1=g_0\circ (\Psi_1\o \mathrm{Id}_\mathfrak{g}\o \mathrm{Id}_\mathfrak{g})+g_0\circ (\mathrm{Id}_\mathfrak{g}\o  \mathrm{Id}_\mathfrak{g} \o\Psi_1)+g_0\circ (\mathrm{Id}_\mathfrak{g}\o \Psi_1 \o \mathrm{Id}_\mathfrak{g})-\Psi_1\circ g_0,\\
&\varphi_{1}-\varphi'_1= \varphi \circ \Psi_1- \Psi_{1}\circ \varphi.
 \end{align*}
That is, we get
  $$((f_1,g_1), \varphi_1)-((f'_1,g'_1), \varphi'_1)=(\delta^1(\Psi_1),-\Phi^1(\Psi_1))=\partial^1(\Psi_1)\in \mathcal{B}_{\mathrm{MDLY}}^{2}(\mathfrak{g}). $$
  Therefore,  $((f'_1,g'_1), \varphi'_1)$ and $((f_1,g_1), \varphi_1)$ are in the same cohomology class in $\mathcal{H}^{2}_{\mathrm{MDLY}}(\mathfrak{g})$.
\end{proof}

\begin{definition}
A  formal deformation $(f_t, g_t, \varphi_t)$  of modified $\lambda$-differential  Lie-Yamaguti algebra  $(\mathfrak{g}, [, ], \{, , \},\varphi)$ is said to be trivial if the deformation $(f_t, g_t, \varphi_t)$ is equivalent
to the undeformed one $(f_0, g_0, \varphi)$.
\end{definition}

\begin{definition}
A modified $\lambda$-differential  Lie-Yamaguti algebra  $(\mathfrak{g}, [, ], \{, , \},\varphi)$ is said to be rigid if every   formal deformation of $\mathfrak{g}$ is trivial deformation.
\end{definition}

\begin{theorem}
 Let $(\mathfrak{g}, [, ], \{, , \},\varphi)$ be a modified $\lambda$-differential  Lie-Yamaguti algebra. If $\mathcal{H}^{2}_{\mathrm{MDLY}}(\mathfrak{g})=0$, then
$(\mathfrak{g}, [, ], \{, , \},\varphi)$ is rigid.
\end{theorem}
\begin{proof}
Let $(f_t, g_t, \varphi_t)$ be a   formal deformation of $(\mathfrak{g}, [, ], \{, , \},\varphi)$. From Proposition \ref{prop:2-cocycle},
$((f_1, g_1), \varphi_1)$ is a 2-cocycle. By $\mathcal{H}^{2}_{\mathrm{MDLY}}(\mathfrak{g})=0$, then there exists a 1-cochain
$$\Psi_1 \in \mathcal{C}^1_{\mathrm{MDLY}}(\mathfrak{g})=\mathcal{C}^1_{\mathrm{LY}}(\mathfrak{g})$$
such that $((f_1,g_1), \varphi_1)= \partial^1(\Psi_1), $  that is, $f_1=\delta_I^1(\Psi_1)$, $g_1=\delta_{II}^1(\Psi_1)$ and $\varphi_1=-\Phi^1(\Psi_1)$.

Setting $\Psi_t = \mathrm{Id}_\mathfrak{g} -\Psi_1t$, we get a deformation  $(\overline{f}_t, \overline{g}_t, \overline{\varphi}_t)$, where
 \begin{align*}
\overline{f}_t=&\Psi_t^{-1}\circ f_t\circ (\Psi_t\otimes \Psi_t),\\
\overline{g}_t=&\Psi_t^{-1}\circ g_t\circ (\Psi_t\otimes \Psi_t\otimes \Psi_t),\\
\overline{\varphi}_t=&\Psi_t^{-1}\circ \varphi_t\circ \Psi_t.
 \end{align*}
 It can be easily verify that $\overline{f}_1=0,\overline{g}_1=0,  \overline{\varphi}_1=0$. Then
    $$\begin{array}{rcl} \overline{f}_t&=& f_0+\overline{f}_2t^2+\cdots,\\
    \overline{g}_t&=& g_0+\overline{g}_2t^2+\cdots,\\
 \varphi_t&=& \varphi_0+\overline{\varphi}_2t^2+\cdots.\end{array}$$
  By Eqs.     \eqref{4.1}--\eqref{4.7},
  we see that $((\overline{f}_2, \overline{g}_2),  \overline{\varphi}_2)$  is still a 2-cocyle, so by induction, we can show
that $(f_t, g_t, \varphi_t)$  is equivalent to the trivial deformation $(f_0, g_0, \varphi)$.  Therefore,  $(\mathfrak{g}, [, ], \{, , \},\varphi)$ is rigid.
\end{proof}

  \section{   Abelian extensions of modified $\lambda$-differential Lie-Yamaguti algebras}\label{sec:extensions}
\def\theequation{\arabic{section}. \arabic{equation}}
\setcounter{equation} {0}
In this section, we consider  abelian extensions of modified $\lambda$-differential Lie-Yamaguti algebras.
We prove that any abelian extension of a modified $\lambda$-differential Lie-Yamaguti algebra  has a representation and a  2-cocycle.  It is further proved that they are classified by the second cohomology.

Let $(\mathfrak{g},\varphi)=(\mathfrak{g}, [, ], \{, , \},\varphi)$  be a modified $\lambda$-differential  Lie-Yamaguti algebra  and $(V, \varphi_V)$ be an abelian modified $\lambda$-differential  Lie-Yamaguti algebra,
 i.e., the Lie-Yamaguti algebra brackets $[,]_V=0$ and $\{,,\}_V=0.$

\begin{definition}
An abelian extension of  $(\mathfrak{g}, \varphi)$ by  $(V, \varphi_V)$
 is  a short exact sequence of   morphisms of modified $\lambda$-differential Lie-Yamaguti algebras
$$\begin{CD}
0@>>> {(V, \varphi_V)} @>i >> (\hat{\mathfrak{g}},[, ]_{\hat{\mathfrak{g}}}, \{, , \}_{\hat{\mathfrak{g}}},\hat{\varphi}) @>p >> (\mathfrak{g}, \varphi) @>>>0
\end{CD}$$
such that  $V$ is an abelian ideal of $\hat{\mathfrak{g}},$  i.e.,  $[u, v]_{\hat{\mathfrak{g}}}=0,\{, u,v\}_{\hat{\mathfrak{g}}}=\{u,v,\}_{\hat{\mathfrak{g}}}=0, \forall u,v\in V$.
\end{definition}

\begin{definition}
 A   section  of an abelian extension $(\hat{\mathfrak{g}},[, ]_{\hat{\mathfrak{g}}}, \{, , \}_{\hat{\mathfrak{g}}},\hat{\varphi})$ of $(\mathfrak{g}, \varphi)$  by  $(V, \varphi_V)$ is a linear map $s:\mathfrak{g}\rightarrow \hat{\mathfrak{g}}$ such that   $p\circ s=\mathrm{Id}_\mathfrak{g}$.
\end{definition}

\begin{definition}
   Let $(\hat{\mathfrak{g}}_1,[, ]_{\hat{\mathfrak{g}}_1}, \{, , \}_{\hat{\mathfrak{g}}_1},\hat{\varphi}_1)$ and  $(\hat{\mathfrak{g}}_2,[, ]_{\hat{\mathfrak{g}}_2}, \{, , \}_{\hat{\mathfrak{g}}_2},\hat{\varphi}_2)$  be two abelian extensions of $(\mathfrak{g},\varphi)$  by  $(V, \varphi_V)$. They are said to be  equivalent if  there is an isomorphism of  modified $\lambda$-differential Lie-Yamaguti algebrass $\eta:(\hat{\mathfrak{g}}_1,[, ]_{\hat{\mathfrak{g}}_1}, \{, , \}_{\hat{\mathfrak{g}}_1},\hat{\varphi}_1)\rightarrow (\hat{\mathfrak{g}}_2,[, ]_{\hat{\mathfrak{g}}_2}, \{, , \}_{\hat{\mathfrak{g}}_2},\hat{\varphi}_2)$
such that the following diagram is  commutative:
\begin{align}
\begin{CD}
0@>>> {(V, \varphi_V)} @>i_1 >> (\hat{\mathfrak{g}}_1,[, ]_{\hat{\mathfrak{g}}_1}, \{, , \}_{\hat{\mathfrak{g}}_1},\hat{\varphi}_1) @>p_1 >> (\mathfrak{g}, \varphi) @>>>0\\
@. @| @V \eta VV @| @.\\
0@>>> {(V, \varphi_V)} @>i_2 >> (\hat{\mathfrak{g}}_2,[, ]_{\hat{\mathfrak{g}}_2}, \{, , \}_{\hat{\mathfrak{g}}_2},\hat{\varphi}_2) @>p_2 >> (\mathfrak{g}, \varphi) @>>>0.\label{5.1}
\end{CD}
\end{align}
\end{definition}

Now for an  abelian extension $(\hat{\mathfrak{g}},[, ]_{\hat{\mathfrak{g}}}, \{, , \}_{\hat{\mathfrak{g}}},\hat{\varphi})$ of $(\mathfrak{g}, \varphi)$  by  $(V, \varphi_V)$  with a section $s:\mathfrak{g}\rightarrow\hat{\mathfrak{g}}$, we define  linear maps $\rho: \mathfrak{g} \rightarrow \mathrm{End}(V)$  and  $\theta, D: \mathfrak{g}\times\mathfrak{g} \rightarrow \mathrm{End}(V)$  by
$$\rho(x)u:=[s(x),u]_{\hat{\mathfrak{g}}},  $$
$$\theta(x,y)u:=\{u,s(x),s(y)\}_{\hat{\mathfrak{g}}},  \quad \forall x,y\in \mathfrak{g}, u\in V.$$
In particular, $D(x,y)u=\theta(y,x)u-\theta(x,y)u=\{s(x),s(y),u\}_{\hat{\mathfrak{g}}}.$
\begin{proposition} \label{prop:representation}
  With the above notations, $(V; \rho, \theta,D, \varphi_V)$ is a representation of the  modified $\lambda$-differential Lie-Yamaguti algebra  $(\mathfrak{g},\varphi)$.
\end{proposition}
\begin{proof}
 In view of \cite{Z15}, $(V; \rho, \theta,D)$ is a representation of  the  Lie-Yamaguti algebra  $\mathfrak{g}$.
 Further, for any $x,y\in \mathfrak{g}$ and $ u\in V,$ ${\hat{\varphi}}s(x)-s(\varphi(x))\in V$ means that $\rho({\hat{\varphi}}s(x))u=\rho(s\varphi(x))u,\theta({\hat{\varphi}}s(x), {\hat{\varphi}}s(y)) u=\theta(s\varphi(x), s(\varphi(y))) u$. Therefore,
we have
 \begin{align*}
&\varphi_V(\rho(x)u)=\varphi_V[s(x),u]_{\hat{\mathfrak{g}}}\\
&=[\varphi_Vs(x),u]_{\hat{\mathfrak{g}}}+[s(x),\varphi_V(u)]_{\hat{\mathfrak{g}}})+\lambda[s(x),u]_{\hat{\mathfrak{g}}}\\
&=\rho(\varphi_V(x))u+\rho(x)\varphi_V(u)+\lambda\rho(x)u,\\
&\varphi_V(\theta(x,  y)u)= \varphi_V\{ u, s(x), s(y)\}_{\hat{\mathfrak{g}}}\\
\quad  &=\{\varphi_V(u), s(x), s(y)\}_{\hat{\mathfrak{g}}}+\{u, \varphi_Vs(x), s(y)\}_{\hat{\mathfrak{g}}}+\{ u, s(x), \varphi_Vs(y)\}_{\hat{\mathfrak{g}}}+2\lambda\{u, s(x), s(y)\}_{\hat{\mathfrak{g}}})\\
\quad  &=\theta(x, y)\varphi_V(u) +\theta(\varphi(x), y)R_Vu+\theta(x, \varphi(y))u+2\lambda\theta(x, y)u.
\end{align*}
Hence,  $(V; \rho, \theta,D, \varphi_V)$ is a representation of   $(\mathfrak{g}, \varphi)$.
\end{proof}

We further define linear maps $\nu:\mathfrak{g}\times \mathfrak{g}\rightarrow V$, $\psi:\mathfrak{g}\times \mathfrak{g}\times \mathfrak{g}\rightarrow V$ and $\chi:\mathfrak{g}\rightarrow V$ respectively by
\begin{align*}
\nu(x,y)&=[s(x), s(y)]_{\hat{\mathfrak{g}}}-s([x, y]),\\
\psi(x,y,z)&=\{s(x), s(y), s(z)\}_{\hat{\mathfrak{g}}}-s(\{x, y,z\}),\\
\chi(a)&=\hat{\varphi}(s(x))-s(\varphi(x)),\quad\forall x,y,z\in \mathfrak{g}.
\end{align*}
We transfer the modified $\lambda$-differential Lie-Yamaguti algebra structure on $\hat{\mathfrak{g}}$ to $\mathfrak{g}\oplus V$ by endowing $\mathfrak{g}\oplus V$ with   multiplications $[, ]_\nu,\{,,\}_\psi$
and a modified $\lambda$-differential
operator  $\varphi_\chi$  defined by
\begin{align}
[x+u, y+v]_\nu&=[x, y]+\rho(x)v-\rho(y)u+\nu(x,y), \label{5.2}\\
 \{x+u, y+v, z+w\}_\psi&=\{x, y, z\}+\theta(y, z)u-\theta(x, z)v+D(x, y)w+\psi(x, y, z),\label{5.3}\\
 \varphi_\chi(x+u)&=\varphi(x)+\chi(x)+\varphi_V(u),  \forall x,y,z\in \mathfrak{g},\,u,v,w\in V.\label{5.4}
\end{align}

\begin{proposition}\label{prop:2-cocycles}
The 4-tuple $(\mathfrak{g}\oplus V,[,]_\nu, \{,,\}_\psi, \varphi_\chi)$ is a modified $\lambda$-differential  Lie-Yamaguti algebra  if and only if
$((\nu,\psi),\chi)$ is a 2-cocycle  of the modified $\lambda$-differential Lie-Yamaguti algebra $(\mathfrak{g}, [, ], \{, , \},\varphi)$ with the coefficient  in $(V; \rho, \theta,D, \varphi_V)$.
 In this case,
$$\begin{CD}
0@>>> {(V, \varphi_V)} @>i >> (\mathfrak{g}\oplus V,[,]_\nu, \{,,\}_\psi, \varphi_\chi) @>p >>(\mathfrak{g},\varphi) @>>>0
\end{CD}$$
 is an abelian extension.
\end{proposition}
\begin{proof}
In view of \cite{Z15},
$(\mathfrak{g}\oplus V,[,]_\nu, \{,,\}_\psi)$ is a  Lie-Yamaguti algebra  if and only if
$\delta^2(\nu,\psi)=0$.

The map $\varphi_\chi$ is a modified modified $\lambda$-differential  operator on $(\mathfrak{g}\oplus V,[,]_\nu, \{,,\}_\psi)$ if and only~if
\begin{align*}
&\varphi_\chi[x + u, y + v]_\nu=[\varphi_\chi(x + u), y + v]_\nu + [x + u, \varphi_\chi(y + v)]_\nu+\lambda[x + u, y + v]_\nu,\\
&\varphi_\chi\{x + u, y + v, z + w\}_\psi=\{\varphi_\chi(x + u), y + v, z + w\}_\psi + \{x + u, \varphi_\chi(y + v), z + w\}_\psi\\
\quad &~~~~~~~~~~~~~~~~~~ ~~~~~~~~~~~~~~~~~~~+\{x + u, y + v,\varphi_\chi(z + w)\}_\psi)+2\lambda\{x + u, y + v, z + w\}_\psi,
\end{align*}
for any $x,y,z\in \mathfrak{g}$ and $u,v,w\in V$. Further,
 we get that the above equations are equivalent to the following equations:
\begin{align}
&\chi[x,y]+\varphi_V(\nu(x,y))=\nu(\varphi(x),y)+\nu(x,\varphi(y))+\lambda\nu(x,y)+\rho(x)\chi(y)-\rho(y)\chi(x), \label{5.5}\\
&\chi\{x,y,z\}+\varphi_V(\psi(x,y,z))=\psi(\varphi(x),y,z)+\psi(x,\varphi(y),z)+\psi(x,y,\varphi(z))+2\lambda\psi(x,y,z)\nonumber\\
&~~~~~~~~~~~~~~~~~~~~~~~~~~~~~~~~~~~~~~~+\theta(y,z)\chi(x)-\theta(x,z)\chi(y)+D(x,y)\chi(z).\label{5.6}
\end{align}
Using Eqs.  \eqref{5.5} and  \eqref{5.6}, we get $\delta_I^{1} (\chi)+ \Phi_I^{2}(\nu)=0$ and $\delta_{II}^{1} (\chi)+ \Phi_{II}^{2}(\psi)=0$  respectively.
Therefore, $ \partial^2((\nu,\psi),\chi)=(\delta^2(\nu,\psi),\delta^1(\chi)+\Phi^{2}(\nu,\psi))=0,$ that is, $((\nu,\psi),\chi)$ is a  2-cocycle.

Conversely, if $((\nu,\psi),\chi)$ is a  2-cocycle  of the modified $\lambda$-differential Lie-Yamaguti algebra $(\mathfrak{g}, [, ], \{, , \},\varphi)$ with the coefficient  in $(V; \rho, \theta,D, \varphi_V)$,  then we have
$\partial^2((\nu,\psi),\chi)=(\delta^2(\nu,\psi),\delta^1(\chi)+\Phi^{2}(\nu,\psi))=0,$ in which $\delta^2(\nu,\psi)=0$, Eqs.  \eqref{5.5} and \eqref{5.6} hold.
So $(\mathfrak{g}\oplus V,[,]_\nu, \{,,\}_\psi, \varphi_\chi)$  is a  modified $\lambda$-differential Lie-Yamaguti algebra.
\end{proof}

\begin{proposition}
Let $(\hat{\mathfrak{g}},[, ]_{\hat{\mathfrak{g}}}, \{, , \}_{\hat{\mathfrak{g}}},\hat{\varphi})$ be an abelian extension of $(\mathfrak{g}, \varphi)$  by  $(V, \varphi_V)$ and $s:\mathfrak{g}\rightarrow\hat{\mathfrak{g}}$   a section. If   $((\nu,\psi),\chi)$ is a 2-cocycle   constructed using the section $s$, then  its cohomology class does not depend on the choice of $s$.
\end{proposition}
\begin{proof}
Let $s_1,s_2:\mathfrak{g}\rightarrow \hat{ \mathfrak{g}}$ be two distinct sections, then we get  two corresponding 2-cocycles $((\nu_1,\psi_1),\chi_1)$ and $((\nu_2,\psi_2),\chi_2)$ respectively. Define
a linear map $\omega: \mathfrak{g}\rightarrow V$ by $\omega(x)=s_1(x)-s_2(x)$ for $x\in\mathfrak{g}$. Further,
\begin{align*}
\nu_1(x, y)=& [s_1(x), s_1(y)]_{\hat{\mathfrak{g}}}-s_1[x, y]\\
\quad= &[s_2(x) + \omega(x), s_2(y) +\omega(y)]_{\hat{\mathfrak{g}}}- s_2([x, y])-\omega[x, y]\\
\quad= & [s_2(x), s_2(y)]_{\hat{\mathfrak{g}}} +\rho(x)\omega(y)-\rho(y)\omega(x)-s_2([x, y])- \omega[x, y]\\
\quad= & \nu_2(x, y) + \delta_{I}^1\omega(x, y),\\
\psi_1(x, y, z)=& \{s_1(x), s_1(y), s_1(z)\}_{\hat{\mathfrak{g}}}-s_1\{ x, y, z\}\\
\quad= &\{s_2(x) + \omega(x), s_2(y) + \omega(y), s_2(z) + \omega(z)\}_{\hat{\mathfrak{g}}}- s_2\{x, y, z\}-\omega\{x, y, z\}\\
\quad=& \{s_2(x), s_2(y), s_2(z)\}_{\hat{\mathfrak{g}}} + \theta(y, z)\omega(x)-\theta(x, z)\omega(y)+D(x, y)\lambda(z)- \omega\{x, y, z\}\\
\quad &-s_2\{x, y, z\}\\
\quad=& \psi_2(x, y, z) + \delta_{II}^1\omega(x, y, z)
\end{align*}
and
\begin{align*}
\chi_1(x) &= \hat{\varphi}s_1(x)- s_1\varphi(x)\\
\quad &= \hat{\varphi}(s_2(x) + \omega(x))-(s_2\varphi(x) +\omega\varphi(x))\\
\quad &= \hat{\varphi}s_2(x)-s_2\varphi(x) + \hat{\varphi}\omega(x)-\omega(\varphi(x))\\
\quad &= \chi_2(x) + \varphi_V \omega(x)-\omega (\varphi(x))\\
\quad &= \chi_2(x)-\Phi^1\omega(x).
\end{align*}
So, $((\nu_1,\psi_1),\chi_1)-((\nu_2,\psi_2),\chi_2)=(\delta^1\omega,-\Phi^1\omega)=\partial^1(\omega)\in  \mathcal{B}_{\mathrm{MDLY}}^{2}(\mathfrak{g},V)$, that is  $((\nu_1,\psi_1),\chi_1)$ and $((\nu_2,\psi_2),\chi_2)$ are in the same cohomological class  in $ \mathcal{H}_{\mathrm{MDLY}}^{2}(\mathfrak{g},V)$.
\end{proof}

\begin{theorem} \label{theorem:classify abelian extensions}
Abelian extensions of a modified $\lambda$-differential Lie-Yamaguti algebra  $(\mathfrak{g}, \varphi)$  by  $(V, \varphi_V)$ are classified by the second cohomology group $\mathcal{H}_{\mathrm{MDLY}}^{2}(\mathfrak{g},V)$ of $(\mathfrak{g}, \varphi)$ with coefficients in the representation $(V; \rho, \theta, D, \varphi_V)$.
\end{theorem}
\begin{proof}
Suppose $(\hat{\mathfrak{g}}_1,[, ]_{\hat{\mathfrak{g}}_1}, \{, , \}_{\hat{\mathfrak{g}}_1},\hat{\varphi}_1)$ and  $(\hat{\mathfrak{g}}_2,[, ]_{\hat{\mathfrak{g}}_2}, \{, , \}_{\hat{\mathfrak{g}}_2},\hat{\varphi}_2)$
 are equivalent abelian extensions   of $(\mathfrak{g}, \varphi)$  by  $(V, \varphi_V)$  with the associated isomorphism $\eta:(\hat{\mathfrak{g}}_1,[, ]_{\hat{\mathfrak{g}}_1}, \{, , \}_{\hat{\mathfrak{g}}_1},\hat{\varphi}_1)\rightarrow (\hat{\mathfrak{g}}_2,[, ]_{\hat{\mathfrak{g}}_2}, \{, , \}_{\hat{\mathfrak{g}}_2},\hat{\varphi}_2)$ such that the diagram in~\eqref{5.1} is commutative.
 Let $s_1$ be a section of $(\hat{\mathfrak{g}}_1,[, ]_{\hat{\mathfrak{g}}_1}, \{, , \}_{\hat{\mathfrak{g}}_1},\hat{\varphi}_1)$. As $p_2\circ \eta=p_1$, we  have
$$p_2\circ(\eta\circ s_1)=(p_2\circ\eta)\circ s_1=p_1\circ s_1= \mathrm{Id}_{\mathfrak{g}}.$$
So the map $\eta\circ s_1$ is a section of $(\hat{\mathfrak{g}}_2,[, ]_{\hat{\mathfrak{g}}_2}, \{, , \}_{\hat{\mathfrak{g}}_2},\hat{\varphi}_2)$. Denote $s_2:=\eta\circ s_1$. Since $\eta$ is an isomorphism of  modified $\lambda$-differential Lie-Yamaguti algebra such that $\eta|_V=\mathrm{Id}_V$, we have
\begin{align*}
\nu_2(x,y)&=[s_2(x), s_2(y)]_{\hat{\mathfrak{g}}_2}-s_2([x,y])\\
&=[\eta(s_1(x)), \eta(s_1(y))]_{\hat{\mathfrak{g}}_2}-\eta(s_1([x, y]))\\
&=\eta\big([s_1(x), s_1(y)]_{\hat{\mathfrak{g}}_1}-s_1([x, y])\big)\\
&=\eta(\nu_1(x,y))\\
&=\nu_1(x,y),\\
\psi_2(x, y, z)&=\{s_2(x), s_2(y), s_2(z)\}_{\hat{\mathfrak{g}}_2}-s_2(\{x, y, z\})\\
\quad &=\eta(\{s_1(x), s_1(y), s_1(z)\}_{\hat{\mathfrak{g}}_1}-s_1\{x, y, z\})\\
\quad &= \psi_1(x, y, z),\\
\chi_2(x)&=\hat{\varphi}(s_2(x))-s_2(\varphi(x))\\
\quad &=\hat{\varphi}(\eta\circ s_1(x))-\eta\circ s_1(\varphi(x))\\
\quad &= \hat{\varphi}(s_1(x))- s_1(\varphi(x))\\
\quad &=\chi_1(x).
\end{align*}
So, all equivalent abelian extensions give rise to the same element in $\mathcal{H}_{\mathrm{MDLY}}^{2}(\mathfrak{g},V)$.

Conversely, given two  cohomologous 2-cocycles $((\nu_1,\psi_1),\chi_1)$ and $((\nu_2,\psi_2),\chi_2)$  in  $\mathcal{H}_{\mathrm{MDLY}}^{2}(\mathfrak{g},V)$,
we can construct two abelian extensions $(\mathfrak{g}\oplus V,[,]_{\nu_1}, \{,,\}_{\psi_1}, \varphi_{\chi_1})$ and  $(\mathfrak{g}\oplus V,[,]_{\nu_2}, \{,,\}_{\psi_2}, \varphi_{\chi_2})$ via Proposition \ref{prop:2-cocycles}. Then  there is  a linear map $\omega: \mathfrak{g}\rightarrow  V$ such that
 $$((\nu_1,\psi_1),\chi_1)-((\nu_2,\psi_2),\chi_2)=\partial^1(\omega)=(\delta^1\omega,-\Phi^1\omega).$$

 Define a linear map $\eta_\omega: \mathfrak{g}\oplus V\rightarrow  \mathfrak{g}\oplus V$ by
$\eta_\omega(x+u):=x+\omega(x)+u$ for $x\in \mathfrak{g}, u\in V.$  Then we have $\eta_\omega$ is an isomorphism of these two abelian extensions.
\end{proof}

\renewcommand{\refname}{REFERENCES}


\begin{thebibliography}{99}



\bibitem{B05} P. Benito, C. Draper, A. Elduque. Lie Yamaguti algebra related to $g_2$. J. Pure Appl. Algebra, 2005,  202, 22--54.

\bibitem{B09} P. Benito, A. Elduque, F. Mart\'{i}n-Herce. Irreducible Lie-Yamaguti algebras.   J. Pure Appl. Algebra, 2009, 213, 795--808.


\bibitem{D22} A. Das. A cohomological study of modified Rota-Baxter algebras. {arXiv Preprint}, 2022, 	arXiv:2207.02273.

\bibitem{G24}  S. Goswami, S. K. Mishra, G. Mukherjee. Automorphisms of extensions of Lie-Yamaguti algebras and inducibility problem. J. Algebra, 2024,   641,  268--306.

\bibitem{Guo23} S. Guo.  Lie-Yamaguti algebras with a derivation.  Acta Math. Sin.(Chin. Ser.), 2023, 66,3, 547--556.



\bibitem{J22}  J. Jiang, Y. Sheng. Cohomologies and deformations of modified $r$-matrices. arXiv Preprint, 2022, arXiv:2206.00411.





\bibitem{K01}  M. K. Kinyon,  A.Weinstein. Leibniz algebras, Courant algebroids and multiplications on reductive homogeneous spaces. Amer. J. Math.,   2001, 123, 525--550.


\bibitem{L15} J. Lin, L. Chen,  Y. Ma.  On the deformaions of Lie-Yamaguti algebras.  Acta. Math. Sin. (Engl. Ser.), 2015, 31, 938--946.


\bibitem{L22} Y. Li, D. Wang.  Cohomology and deformation theory of modified Rota-Baxter Leibniz algebras. arXiv Preprint, 2022,	arXiv:2211.09991.

\bibitem{L23} Lin  J.,  Chen  L.  Quasi-derivations of Lie-Yamaguti algebras. J. Alg. Appl., 2023, 22(5), 2350119.

\bibitem{M22} B. Mondal, R. Saha.   Cohomology of modified Rota-Baxter Leibniz algebra of weight $\kappa$.  {arXiv Preprint}, 2022,	arXiv:2211.07944.


\bibitem{N54}  K.  Nomizu.  Invariant affine connections on homogeneous spaces. Amer. J. Math., 1954, 76, 33--65.

\bibitem{Peng}  X. Peng, Y. Zhang, X. Gao, Y. Luo. Universal enveloping of (modified) $\lambda$-differential Lie algebras.  Linear Mult. Alg., 2022, 70, 1102--1127.



\bibitem{S21}  Y. Sheng,  J. Zhao,  Y. Zhou.  Nijenhuis operators, product structures and complex structures on Lie-Yamaguti algebras.   J. Alg. Appl.,  2021,  20(8),  2150146.


\bibitem{S22}  Y. Sheng,  J. Zhao.  Relative Rota-Baxter operators and symplectic structures on Lie-Yamaguti algebras. Commun. Algebra, 2022, 50(9), 1--18.

\bibitem{Semonov}  M. Semenov-Tian-Shansky.   What is a classical r-matrix?  Funct. Anal. Appl., 1983, 17, 259--272.


\bibitem{Teng23}  W. Teng, F. Long, Y. Zhang. Cohomologies of modified $\lambda$-differential Lie triple systems and applications. AIMS Math., 2023, 8(10), 25079--25096.

\bibitem{Teng23-2}  W. Teng, J. Jin, F. Long. Generalized reynolds operators on Lie-Yamaguti Algebras.   Axioms, 2023, 12(10),  934.

\bibitem{Teng23-3}  W. Teng, H. Zhang. Deformations and extensions of modified $\lambda$-differential 3-Lie algebras. Mathematics, 2023, 11(18), 3853.

\bibitem{Teng24}  Y. Xiao, W. Teng. Representations and cohomologies of modified $\lambda$-differential Hom-Lie algebras. AIMS Math., 2024, 9(2),  4309--4325.

\bibitem{TG24}  W. Teng, S. Guo. Modified Rota-Baxter Lie-Yamaguti algebras.  {arXiv Preprint}, 2024, 	arXiv:2401.17726.  DOI: 10.13140/RG.2.2.12760.67845


\bibitem{Y57} K. Yamaguti.  On the Lie triple system and its generalization.  J. Sci. Hiroshima Univ. Ser. A, 1957/1958,  21, 155-160.

\bibitem{Y67}  K. Yamaguti.  On cohmology groups of general Lie triple systems.  Kumamoto J. Sci. A,  1967/1969,  8,  135-146.

\bibitem{Z22}  J. Zhao, Y. Qiao. Cohomologies and deformations of relative Rota-Baxter operators on Lie-Yamaguti algebras. Mathematics,  2024, 12(1),  166

\bibitem{Z15}  T. Zhang,  J. Li.  Deformations and extensions of Lie-Yamaguti algebras.  Linear Mult. Alg., 2015, 63, 2212-2231.


































\end{thebibliography}
\end{document}